\documentstyle[12pt]{article}
\newlength{\abstractwidth}
\renewcommand{\thefootnote}{\fnsymbol{footnote}}
\renewcommand{\thanks}[1]{\footnote{#1}} 
\newcommand{\starttext}{ \setcounter{footnote}{0}
\renewcommand{\thefootnote}{\arabic{footnote}}}

\newcommand{\be}{\begin{equation}}
\newcommand{\bea}{\begin{eqnarray}}
\newcommand{\eea}{\end{eqnarray}} \newcommand{\ee}{\end{equation}}
 \newcommand{\<}{\langle}
\renewcommand{\>}{\rangle} \def\ba{\begin{eqnarray}}
\def\ea{\end{eqnarray}}

\def\o{\omega}
\def\Re{{\rm Re}}

\def\log{\,{\rm log}\,}

\def\o{\omega}

\def\e{\varepsilon}

\def\o{\omega}

\def\na{\nabla}

\def\p{\partial}

\def\ddb{{\partial\bar\partial}}

\def\na{{\nabla}}

\def\[{{\bf [}}
\def\]{{\bf ]}}

\begin{document}
\starttext \baselineskip=18pt \setcounter{footnote}{0}
\newtheorem{theorem}{Theorem}
\newtheorem{lemma}{Lemma}
\newtheorem{corollary}{Corollary}
\newtheorem{definition}{Definition}
\newtheorem{conjecture}{Conjecture}
\newtheorem{proposition}{Proposition}

\begin{center}
{\Large \bf GEOMETRIC FLOWS FROM UNIFIED STRING THEORIES
\footnote{Contribution to Surveys in Differential Geometry, Vol. 27 (2022), ``Forty Years of Ricci Flow", edited by H.D.  Cao, R. Hamilton, and S.T. Yau. Work supported in part by the National Science Foundation under grants DMS-18-55947 and DMS-22-03273.}}

\medskip
\centerline{Duong H. Phong}

\medskip

\begin{abstract}

A survey of new geometric flows motivated by string theories is provided. Their settings can range from complex geometry to almost-complex geometry to symplectic geometry. From the PDE viewpoint, many of them can be viewed as intermediate flows between the Ricci flow and the K\"ahler-Ricci flow, albeit often coupled to flows of additional fields. In particular, a survey is given of joint works of the author with Tristan Collins, Teng Fei, Bin Guo, Sebastien Picard, and Xiangwen Zhang.

\end{abstract}

\end{center}

\baselineskip=15pt
\setcounter{equation}{0}
\setcounter{footnote}{0}

\section{Introduction}
\setcounter{equation}{0}

The four known fundamental forces of nature are each described by either the Einstein equation, or by the Yang-Mills equation, both of which have had a profound influence on the development of geometry and the theory of partial differential equations. A grand dream of theoretical physics is however a unified theory of all interactions, for which string theories, unified themselves with $11$-dimensional supergravity since the mid 1990's into M Theory, remain at this moment as the only viable candidate. Presumably M Theory will be governed by equations which unify in a suitable sense the Einstein equation and the Yang-Mills equation, and which will surely be of great importance for mathematics. The problem is that
M Theory is still far from being well understood, and is even at this moment only known indirectly through its various limits. 

\medskip
As a first step, we can still  consider the equations for each of the limiting theories of M Theory, and in fact, of their low-energy limits, which are $10d$ or $11d$ theories of point particles. Historically, this was actually done first, with the identification in 1985 by Candelas, Horowitz, Strominger, and Witten \cite{CHSW} of Calabi-Yau 3-folds as vacua for  supersymmetric compactifications of the heterotic string to a $4$-dimensional space-time. At the time, Calabi-Yau 3-folds appeared so rigid as to preclude searches for other solutions of the heterotic string or for the other string theories. But the situation has changed a lot since, especially with dualities and non-perturbative effects, and it is imperative to try and understand these other vacua. There is now strong evidence that each of these theories leads to a new notion of {\it special geometry}, and is governed by new partial differential equations which are of considerable interest in their own right. Thus the above Calabi-Yau 3-folds should be interpreted as manifestations of $SU(3)$ holonomy with respect to the Chern unitary connection, but this may be only one among many special geometries yet to be discovered.

\medskip
Remarkably, geometric flows turn out to be the most effective method to investigate these new equations. An early example is the $G_2$ flow of Bryant \cite{Bry, BryXu}, which can be viewed as an odd-dimensional version of the K\"ahler-Ricci flow. More generally, it turns out that the equations for the limiting theories of M Theory all involve a cohomological constraint. In the absence of a $\p\bar\p$-Lemma, the most effective way of implementing a constraint is to start from an initial data satisfying it, and let the data evolve by a geometric flow preserving the constraint.

\medskip
In this paper, we shall provide a survey of the geometric flows which have appeared recently in unified string theories, focusing on the PDE aspects. The original formulation of these flows varies a great deal: some flows arise as flows of $(2,2)$-forms on a 3-dimensional complex manifold (\S 2 and \S 4), others as flows of real $3$-forms on a symplectic 6-manifold (\S 3), others as yet as parabolic reductions of $11d$ supergravity (\S 5), or as flows of spinors (\S 6). But they all induce flows of the Riemannian structure $g_{jk}$, and we can try and compare them from there. Typically they are then of the form
\bea
\p_t g_{jk}=-2 e^u R_{jk}+\cdots
\eea
and involve couplings with the flows of other fields. Here $\cdots$ represent some lower order terms, which may also involve reparametrizations and/or gauge transformations. The appearance of the Ricci tensor is not surprising in a reparametrization invariant flow of second order. Thus the distinctive features of the flow are rather to be found in the lower order terms and/or in the additional fields to which the metric $g_{jk}$ is coupled. It should also be kept in mind that, while cohomological constraints and any underlying special geometry may not be evident from the above flows of the metrics $g_{jk}$, they are there and of considerable importance, even though they are usually much weaker than in K\"ahler geometry. In this sense, the geometric flows of unified string theories can be viewed as intermediates between the Ricci flow and the K\"ahler-Ricci flow.

\medskip
Geometric flows from unified string theories have begun to be investigated systematically only a few years ago, and very little is as yet known about them.
While we can build on the many powerful techniques and insights developed over the years  for the study of classic flows such as the Ricci flow \cite{Ha82, Ha95} and the Yang-Mills flow \cite{UY, D}, the more recent flows invariably lead to new difficulties of their own. In all likelihood, they will require new methods which may well prove to be valuable in their own right for the theory of partial differential equations. 
There is much then to understand and discover, and this whole area should be very fertile for research. With this in mind, we have thrived in this survey to present the new flows in as mathematically self-contained a way as possible, so that the readers can begin studying them as pure PDE questions if they choose. The physics context for these equations is especially rich and complex, and we have to content ourselves with providing some references to the literature, which is immense.

\section{The Type IIB flow}
\setcounter{equation}{0}

We begin with some general remarks, which apply not just to this section, but to the other sections as well. String theories and M Theory \cite{BBS, HW, W, To} are theories of extended objects, but for our purposes, we only consider limits where they reduce to theories of point particles, in dimensions 10 and 11 respectively. These limits are supergravity theories, in the sense that they incorporate gravity and are supersymmetric, i.e., they admit a symmetry exchanging tensor fields (``bosons") and spinor fields (``fermions"). Gravity is described by a metric which is a symmetric rank-two tensor, and its supersymmetric partner is described by a field called the gravitino field, which is a one-form valued in the bundle of spinors. We are interested in supersymmetric vacua of these theories, that is, space-times whose fields satisfy the field equation, and whose gravitino field remains invariant under supersymmetry transformations. This last requirement reduces to the existence of a spinor field which is covariantly constant under a connection with possibly non-vanishing flux (see \S 7 for some more details). It is one of the key new requirements which set string theories and M Theory apart from previous gauge theories and gravity theories. 

\smallskip
Except for $11d$ supergravity, all the equations considered in the present survey can be traced back to the requirement of existence of a covariantly constant spinor. This implies in particular a reduced holonomy. For phenomenological reasons, it is desirable to compactify space-time to $M^{3,1}\times X$, where $M^{3,1}$ is Minkowski or a maximally symmetric $4$-dimensional space-time, and $X$ is an internal space of dimension either $6$ or $7$, depending on whether we compactify string theories or $11$-dimensional supergravity. It is also desirable
to preserve supersymmetry upon compactification. Thus our considerations can apply to the original space-time or more often, the internal space $X$.

\smallskip
The equations for string theories have been extensively explored in the literature
\cite{Gr,T,GMW,GGHPR,GP, GMR, Cardoso, Gurrieri}, notably from the point of view of reduced holonomy and G structures. In the case of a six-dimensional internal space $X$, we are then in the case of $SU(3)$ structures, for which a detailed classification of what the flux can be \cite{Fr, CS}, and this has proved to be a powerful tool for a classification of possible supersymmetric vacua.

\medskip
Our goal is to identify the equations which exhibit the typical features of each of the string theories. For this, we rely on the work of Tseng and Yau \cite{TY1, TY2, TY3}, who have proposed some representative models.
For the Type IIB string (compactified to a six-dimensional internal space $X$), the setting is then as follows.

\medskip
Let $X$ be a compact $3$-dimensional complex manifold, equipped with a nowhere vanishing holomorphic $(3,0)$-form $\Omega$. A supersymmetric compactification of the Type IIB string with an O5/D5 brane source can be described by a Hermitian metric $\o$
satisfying the following system of equations
\bea
d\omega^2=0, \qquad i\partial\bar\partial (|\Omega|_\omega^{-2}\omega)=\rho_B
\eea
where $|\Omega|_\omega$ is defined by $i\Omega\wedge\bar\Omega=|\Omega|_\omega^2\omega^3$, and $\rho_B$ is the Poincare dual of a linear combination of holomorphic $2$-cycles. 

\smallskip
We observe that the equation $d\o^2=0$ is reminiscent of a K\"ahler condition, but it is a lot weaker. Metrics satisfying this condition are sometimes called ``balanced" in the mathematics literature \cite{Michelsohn}. In particular, there is no known analogue of the $\p\bar\p$ Lemma which can provide an effective parametrization of all forms $\o$ with $\o^2$ closed. In the absence of such a parametrization, we try and implement the balanced condition by introducing instead the following flow of $(2,2)$-forms
\bea
\p_t(\o^2)=i\partial\bar\partial (|\Omega|_\omega^{-2}\omega)-\rho_B
\eea
with any initial data $\o_0$ satisfying $d\o_0^2=0$. The point is that the right-hand side of the flow is a closed form, so the closedness of the form $\o^2$ is preserved for all time, and formally, if the flow exists for all time and converges, then its stationary points provide solutions to the full Type IIB system.

\smallskip
It is shown in \cite{PPZ1} that the leading terms in this equation satisfy Hamilton's conditions for the short-time existence of the flow for smooth data.  We also observe that, in this form, the Type IIB flow admits a ready generalization to any compact complex manifold of dimension $m$,
\bea
\p_t(\o^{m-1})
=
i\partial\bar\partial (|\Omega|_\omega^{x}\omega^{m-2})-\rho_B
\eea
which may be of interest in its own right. Here $\rho_B$ is now the Poincare dual of a linear combination of holomorphic $m-1$-cycles. It is also not difficult to rewrite this Type IIB flow as a flow of $(1,1)$-forms \cite{PPZ2}.

\subsection{The case $\rho_B=0$}

We describe now some precise results about the above flow in the case of no source $\rho_B=0$. In this case, even for general complex dimension $m$, there is a remarkable change of variables that allows the dependence of the Type IIB flow on the form $\Omega$ to be completely eliminated from the flow, and reside only in the initial data \cite{FP}. Indeed set 
\bea
\eta=|\Omega|_\o\o
\eea
Then it is shown in \cite{FP} that the Type IIB flow can be expressed in terms of $\eta$ as
\bea
i^{-1}\p_t\eta_{\bar kj}={1\over m-1}(-R_{\bar kj}+{1\over 2}T_{j\bar m p}\bar T_{\bar k}{}^{\bar mp})
\eea
where $R_{|bar kj}$ is the Ricci-Chen tensor, for initial data $\eta_0$ satisfying the conformally balanced condition $d(|\Omega|_\eta\eta^{m-1})=0$. The above right hand side actually coincides with a specific flow within the family of flows introduced in \cite{ST}, as a natural generalization in Hermitian geometry of the K\"ahler-Ricci flow. Remarkably, it is the specific flow identified by Ustinovskyi \cite{U} as preserving the Griffiths and Nakano-positivity of the tangent bundle. We note however that, in our context, it is essential that the initial data be conformally balanced. The following Shi-type estimates were established in \cite{FPPZb}:

\begin{theorem}
\label{thm:TypeIIBShi}
{\rm \cite{FPPZb}}
Assume that the Type IIB flow exists on the time interval $[0,T)$, and that
\bea
K_1^{-1}\omega(0)\leq \omega(t) \leq K_1\omega_0,
\qquad
|d\log |\Omega||\leq K_2
\eea
Then there exists a constant $C$ and $0<\alpha<1$, depending on $K_1,K_2$, $\Psi$, and the initial data, but not on $T$, so that
\bea
\|g\|_{C^{2+\alpha,1+{\alpha\over 2}}(X\times [0,T))}\leq C
\eea
In particular, the flow extends then across $T$.
\end{theorem}

Traditionally, the $C^3$ estimates for the Monge-Amp\`ere equation and the K\"ahler-Ricci flow relied on a complex version of an identity going back in Calabi. Here, no potential for a K\"ahler metric is available, but we can use rely instead on the method of \cite{PSS}, which gives estimates for the connection $\na h h^{-1}$, where $h$ is the relative endomorphism between a reference metric and the evolving metric. Applications of the above theorem to the formation of singularities for the Type IIB flow can be found in \cite{K1}.

\smallskip

The Type IIB flow should provide a powerful probe for the geometry of the underlying complex manifold $X$. As a first test, we can show that it is certainly capable of giving a new proof, with new estimates, of Yau's fundamental theorem \cite{Y} on the existence of Ricci-flat K\"ahler metrics on K\"ahler manifolds with $c_1(X)=0$. 
Thus assume that $\hat\chi$ is a K\"ahler metric on $X$ and take as initial data for the Type IIB flow a metric $\omega(0)$
satisfying $| \Omega |_{\omega(0)} \omega(0)^{n-1} = \hat{\chi}^{n-1}$. Then it can be readily shown \cite{PPZ4} that
the Type IIB flow exists for all time $t>0$, and as $t \rightarrow \infty$, the solution $\omega(t)$ converges smoothly to a K\"ahler, Ricci-flat, metric $\omega_\infty$.

\smallskip

The point in this case is that the $(1,1)$-form $\chi$ defined by
\bea
| \Omega |_{\omega(t)} \omega(t)^{n-1} = \chi(t)^{n-1}
\eea
will remain K\"ahler along the flow. If we set $\chi(t) = \hat{\chi} + i \ddb \varphi(t)$, 
we can rewrite the Type IIB flow as a parabolic Monge-Amp\`ere flow in $\varphi$,
\bea
\p_t \varphi = e^{-f} {(\hat\chi+i\partial\bar\partial\varphi)^n
\over\hat\chi^n}, \ \ \varphi(x,0)=0
\eea
for a suitable given function $f$, subject to the plurisubharmonicity condition $\hat{\chi} + i \ddb \varphi >0$.

\medskip

Note that, unlike the K\"ahler-Ricci flow or the inverse Monge-Amp\`ere flow \cite{CHT, CaK}, this flow is not concave in $D^2\varphi$. Because of this, we need a new way of obtaining $C^2$ estimates. It turns out that a new test function, different from the classical one used by Yau and Aubin for the $C^2$ estimates for the Monge-Amp\`ere equation,
{$$
G(z,t)=\log {\rm Tr}\, h-A(\varphi-{1\over[\hat\chi^n]}\int_X\varphi\hat\chi^n)+B[{(\hat\chi+i\partial\bar\partial\varphi)^n\over\hat\chi^n}]^2
$$}
is needed. The new additional term is the square of the Monge-Amp\`ere determinant.
It has proved to be useful in $C^2$ estimates for other non-linear parabolic flows as well \cite{PZ, KJS}, and confirms that tools developed for the flows of unified string theories will also be of interest for other PDE's.

\medskip
It is not difficult to see that the only stationary points of the Type IIB flow with no source are Ricci-flat K\"ahler metrics. Thus the flow can help determine whether the conformally balanced manifold $X$ is actually K\"ahler. For this, the following stability theorem due to Bedulli and Vezzoni \cite{BVb} is of particular interest:
\begin{theorem}
\label{thm:BV} {\rm \cite{BVb}}
If $(X,\Omega,\chi)$ is a K\"ahler manifold with $\chi$ a K\"ahler Ricci-flat metric, then for any $\epsilon>0$, there exists $\delta>0$ so that, if the initial data $\omega_0$ satisfies
\bea
| |\Omega|_{\omega_0}\omega_0^{n-1}-\chi^{n-1}|<\delta
\eea
then the Type IIB flow exists for all time,  $| |\Omega|_{\omega(t)}\omega(t)^{n-1}-\chi^{n-1}|<\epsilon$, and 
$|\Omega|_\omega\omega^{n-1}$ converges to $\omega_\infty^{n-1}$ with $\omega_\infty$ an astheno-Kaehler metric. If $\omega_0$ is conformally balanced, then $\omega_\infty$ is K\"ahler Ricci-flat.
\end{theorem}

Note that this stability theorem allows in particular to consider more general initial data than conformally K\"ahler.
Its proof is based on a general method of the authors \cite{BV, BVa} for proving the stability of flows which exist by the Hamilton-Nash-Moser implicit function theorem.

\subsection{The case with source $\rho_B$}

We discuss briefly the case of non-vanishing sources $\rho_B$. If the source $\rho_B$ is a smooth closed $(2,2)$-form, then the short-time existence of the flow follows from the same arguments as in the case of no source, and Theorem 1 still applies. 

\smallskip

The more interesting case is for the source to have singularities. In general, the theory of geometric flows and canonical metrics is much less developed in the cases with singularities than as in the smooth case. For example, only very recently has the existence of Hermitian-Einstein metrics been established on normal K\"ahler spaces in the sense of Grauert, and this required new analytic tools \cite{Caoetal}. In the case of the Type IIB flow, practically nothing is known in the case with singular $\rho_B$. If $\rho_B$ is given, we may naively try and consider the Type IIB flow with a source given by a regularization of $\rho_B$, and determine the effects of removing the regularization. But if $\rho_B$ has Dirac-type singularities, we may hope that such Dirac-type singularities would arise dynamically from the flow, when it does not admit smooth critical points. As we shall see in the discussion of the Type IIA flow with source $\rho_A$ below, there does not appear to be any other viable for such a source $\rho_A$ to develop.

\section{The Type IIA flow}
\setcounter{equation}{0}

Again, we focus on a representative case of supersymmetric compactfification of the Type IIA string with brane sources, as proposed by Tseng and Yau \cite{TY1, TY2, TY3}.

\smallskip

Let $X$ be this time a real $6$-dimensional compact symplectic manifold, in the sense that it admits a closed, non-degenerate $2$-form $\omega$ (but there may be no compatible complex structure, so it may not be a K\"ahler form). Given a real $3$-form $\varphi$, let $J_\varphi$ be the almost-complex structure associated some 20 years ago to $\varphi$ by Hitchin \cite{Hitb}. Then the Type IIA equation is the following system of equations for
a real positive and primitive $3$-form $\varphi$,  
\bea
dd^\Lambda(|\varphi|^2\star\varphi)=\rho_A,
\quad d\varphi=0.
\eea
Here $\varphi$ is positive means that the quadratic form
$g_\varphi(X,X)=\omega(X,J_\varphi X)$ is a metric, $\rho_A$ is the Poincare dual of a linear combination of special Lagrangians, $d^\Lambda=d\Lambda-\Lambda d$ is the symplectic adjoint of $d$, and $\star$ is the Hodge star operator of the metric $g_\varphi$.

\medskip
As in the Type IIB flow, a seemingly innocuous condition, but in practice rather difficult to implement, is $d\varphi=0$. And just as in the previous case of the Type IIB flow, this suggests looking for solutions of the above Type IIA system as stationary points of the following flow
\bea
\label{TypeIIA flow}
\partial_t\varphi=dd^\Lambda(|\varphi|^2\star\varphi)-\rho_A,
\eea
with an initial value $\varphi_0$ which is closed and primitive. Since the right hand side is closed and primitive, formally the closedness and primitiveness of the initial condition will be preserved, and the stationary points of the flow are automatically the solutions of the Type IIA system. 

\subsection{Type IIA geometry}

Hitchin's construction of an almost-complex structure $J_\varphi$ from any non-degenerate $3$-form $\varphi$ on any $6$-dimensional manifold $M$, was a pointwise and purely algebraic construction. The above setting for the Type IIA string requires several specific additional properties: the form $\varphi$ has to be closed, the manifold $M$ is also equipped with a symplectic structure, so it makes sense to require that $\varphi$ be primitive. These few specific requirements actually lead to a very rich geometric structure, which is called Type IIA geometry, and whose key features are the following.

\smallskip
Set $\Omega_\varphi=\varphi+iJ_\varphi\varphi$. Then we have
\bea
{\cal D}^{0,1}\Omega_{\varphi}=0
\eea
so that $\Omega_\varphi$ is formally a holomorphic form with respect to the almost-complex structure $J_\varphi$. Next the Nijenhuis tensor has only $6$ independent components, and in this sense, the almost-complex structure is not far from integrable. But the most important property of Type IIA geometries is that they have $SU(3)$ holonomy with respect to a specific connection, namely the projected Levi-Civita connection $\tilde{\cal D}$ with respect to the almost-complex structure $J_\varphi$ and the metric 
\bea
\label{eqn:rescaled}
\tilde g_\varphi=|\varphi|^2g_\varphi.
\eea
More specifically, recall that given an almost-complex structure $J$ and Nijenhuis tensor $N$, then any metric $g(X,Y)$ satisfying the compatibility condition $g(JU,JV)=g(U,V)$ determines a 2-form $\o(U,V)=g(U,JV)$ and a line of unitary connections ${\cal D}^t$, called the Gauduchon line \cite{Gau}, preserving $J$ and $g(X,Y)$, given explicitly by
\bea
{\cal D}_i^tX^m=
\na_iX^m+g^{mk}(-N_{ijk}-V_{ijk}+tU_{ijk})X^j
\eea
where $\na$ is the Levi-Civita connection of the metric $g(X,Y)$, and the tensors $U_{ijk}$ and $V_{ijk}$ are defined by
\bea
U^m{}_{bc}
={1\over 4}((d^c\o)^m{}_{bc}+
(d^c\o)^m{}_{jk}J^j{}_bJ^k{}_c),
\quad
V^m{}_{bc}
={1\over 4}((d^c\o)^m{}_{bc}-
(d^c\o)^m{}_{jk}J^j{}_bJ^k{}_c).
\eea
Here $d^c\o=-Jd\o$, and is given explicitly in components by
$(d^c\o)_{abc}=-J^k{}_cJ^j{}_bJ^i{}_a(\p_i\o_{jk}+\p_j\o_{ki}+\p_i\o_{jk})$. Note that all three tensors $N,U,V$ are anti-symmetric in the last two indices. 
When $d\o=0$, the Gauduchon line reduces to a single connection. But in general, 
there are several connections which are of noteworthy: $t=1$ corresponds to the Chern connection, $t=-1$ is the Bismut connection, characterized by the property that its torsion is a $3$-form. and $t=0$ defines the so-called {\it projected Levi-Civita connection}.

\smallskip
We return now to our discussion of Type IIA geometries on a symplectic $6$-dimensional manifold $M$. Their most basic property is given in the following theorem \cite{FPPZ}:

\begin{theorem}
\label{thm:SU3}
{\rm\cite{FPPZ}}
Let $(J_\varphi,\varphi, \o)$ be a Type IIA geometry, and let $\tilde g_\varphi$ be the rescaled metric defined by (\ref{eqn:rescaled}). Let $\tilde\o=|\varphi|^2\o$ be the corresponding rescaled $2$-form. Then we have
\bea
{\cal D}^0({\Omega_\varphi\over|\Omega_\varphi|})=0
\eea
where ${\cal D}^0$ is the projected Levi-Civita connection of the metric $\tilde g_\varphi$, with respect to the almost-complex structure $J_\varphi$. Thus the Type IIA geometries have SU(3) holonomy, but with respect to the connection ${\cal D}^0$. 
\end{theorem}

We observe that the original symplectic form $\o$ is closed, so the Gauduchon line with respect to the metric $g_\varphi$ reduces to a single connection. However, because of the rescaling, the form $\tilde\o$ is not closed, and we do have a whole line of unitary connections with respect to the metric $\tilde g_\varphi$. The above theorem singles out the projected Levi-Civita connection as the natural connection for Type IIA geometries.

\subsection{The case of $\rho_A=0$}

The case of no source is already geometrically interesting and presents many new analytic challenges. First, we have to establish the short-time existence of the flow:

\begin{theorem}
\label{thm:short}
{\rm \cite{FPPZ}}
The Type IIA flow with $\rho_A=0$ always exists for at least a short-time for any initial data $\varphi$ which is closed, positive, and primitive. The properties of closedness and primitiveness are preserved along the flow.
\end{theorem}

The new difficulty here is that, upon reparametrization by a vector field $V$ to lift the degeneracy inherent to the weak parabolicity of the Type IIA flow, the symplectic form $\o$ also evolves in time. Thus we obtain a coupled flow of $(\varphi,\o)$ which remains only weakly parabolic. The remedy turns out to introduce the following  {\it regularized} coupled flow of $(\varphi,\o)$,
\bea
\partial_t\varphi
=d\Lambda d(|\varphi|^2J\varphi)-BdJd(|\varphi|^2\Lambda J\varphi)+d(\iota_V\varphi),
\quad
\p_t\omega=d(\iota_V\omega)
\eea
for a fixed positive constant $B$. Applying Hamilton's Nash-Moser arguments, this flow can be shown to admit short-time existence. Furthermore, it can be shown to preserve the primitiveness of $\varphi$. But for primitive data, the regularization term with coefficient $B$ vanishes, and thus the solution of the reparametrized, regularized flow satisfies the original Type IIA flow.

\smallskip
Further analysis of the Type IIA flow requires working out the corresponding flow of metrics:

\begin{theorem}
\label{thm:flowofmetrics}
{\rm \cite{FPPZ}}
Assume that $\varphi$ flows by the Type IIA flow. There are two ways of expressing the flow of metrics induced by the Type IIA flow. 

{\rm (a)} Set $\check g_\varphi=|\varphi|^{-2}g_\varphi$ and $u=\log|\varphi|^2$ .  Then
\bea
\partial_t(\check g_\varphi)_{ij}
&=&
e^{{3\over 2}u}\big\{-2\check R_{ij}+{3\over 2}u_iu_j-u_{Ji}u_{Jj}+4u^k(N_{ikj}+N_{jki})-4N^{kp}{}_iN_{pkj}
\nonumber\\
&&+{1\over 2}(|du|^2_{\check g}+|N|_{\check g}^2)\check g_{ij}\big\}
\eea
Note that the volume of $g_\varphi$ is known and equal to $\o^3/3!$, so $u$ is determined by $\check g_{\varphi}$, and the above flow is self-contained.

{\rm (b)} Alternatively, we can consider the flow of the pair $(g_\varphi,u)$. If $\varphi$ flows by the Type IIA flow, then the pair flows by
\bea
\p_tg&=&e^u[-2R_{ij}+2\na_i\na_ju-4(N^2_-)_{ij}
+u_iu_j-u_au_bJ^a{}_iJ^b{}_j+
4u_p(N_i{}^p{}_j+N_j{}^p{}_i)]\nonumber\\
\p_tu&=&e^u\Delta u+e^u(2|\na u|^2+|N|^2)
\eea
\end{theorem}

\smallskip
An important feature of the Type IIA flow which distinguishes it from many known flows of almost-complex structures, e.g. \cite{LW}, is that the norm of the Nijenhuis tensor in the Type IIA flow itself obeys a well-behaved flow,
\bea
(\p_t-e^u\na)|N|^2
&=&e^u[-2|\na N|^2+(\na^2u)\star N^2+Rm\star N^2\nonumber\\
&&
+N\star\na N\star(N+\na u)+|N|^4+N^2\star \na u+N^2\star(\na u)^2].
\eea

\medskip
From here, we can establish Shi-type estimates which show that the singularities of the Type IIA flow can be traced to only a finite number of geometric quantities blowing up:

\begin{theorem}
\label{thm:Shi}
{\rm \cite{FPPZ}}
Assume that the Type IIA flow exists on some time interval $[0,T)$ and that
\bea
|\log \varphi|+|Rm|\leq A
\eea
for some constant $A$.
Then for any $\alpha$, we have
\bea
|\nabla^\alpha\varphi|\leq C(A,\alpha, T,\varphi(0))
\eea
In particular, if $[0,T)$ is the maximum time internal of existence, we must have
\bea
{\rm lim}_{t\to T}{\sup}_X(|\log\varphi|+|Rm|)=\infty
\eea
\end{theorem}

\medskip
Singularity models resulting from these estimates are considered in \cite{K2}.
The stationary points of the Type IIA flow correspond to integrable complex structures, and the resulting metrics are Ricci-flat K\"ahler metrics. 
As shown in \cite{FPPZb}, the Type IIA flow is dynamically stable in the following sense:

\begin{theorem}
{\rm \cite{FPPZb}}
Let $(\bar\varphi,\bar\o)$ be a Ricci-flat Type IIA structure on a compact oriented $6$ manifold $M$. Then there exists $\e_0>0$ with the following property. For any Type IIA structure $(\varphi_0,\o_0)$ satisfying
\bea
|\o_0-\bar\o|_{W^{10,2}}+|\varphi_0-\bar\varphi|_{W^{10,2}}<\e_0,
\eea
the Type IIA flow with initial value $(\varphi_0,\o_0)$ exists for all time, and converges in $C^\infty$ to a Ricci-flat Type IIA structure $(\varphi_\infty,\o_0)$.
\end{theorem}

As an immediate consequence, if a compact symplectic $6d$ manifold $(M,\bar\o)$ admits a compatible integrable complex structure, then any sufficiently close symplectic form will also admit a compatible integrable complex structure.

\smallskip

More generally, it can be seen in explicit examples that the Type IIA flow will evolve towards almost-complex structures which are optimal in a suitable sense. 
First we consider the Type IIA flow on the nilmanifolds constructed by de Bartholomeis and Tomassini \cite{deBT}. There the basic structure is a nilpotent Lie group, and the natural ansatze for $\varphi$ are preserved, and reduce the Type IIA flow as a system of ODE's which can be solved explicitly. We find in this way examples where the flow exists for all times, but the limit {${\rm lim}_{t\to \infty}J_t$} does not exist, but
\bea
{\rm lim}_{t\to\infty}|N|^2=0.
\eea

Next, another instructive example is provided by the symplectic half-flat structures on the solvmanifold introduced by Tomassini and Vezzoni \cite{TV}. Some natural ansatze for $\varphi$ are again preserved by the Type IIA flow, which reduces then to a system of ODE's which can be solved exactly. We find then the following interesting phenomena:

\smallskip

$\bullet$ For any initial data, the flow for $\varphi$ develops singularities in finite time. However, the limits of $J_\varphi$ and $g_\varphi$ continue to exist. 

$\bullet$ For certain initial data, $\varphi$ is a self-expander, while $J$ is stationary, and in fact a critical point of the energy functional of Blair-Ianus \cite{BI} and Le-Wang \cite{LW}.

$\bullet$ For other initial data, as $t$ approaches the maximum time of existence $T$, the limit {${\rm lim}_{t\to T}J$} exists, and is a harmonic almost-complex structure, and in particular a minimizer for $|N|^2$. See also \cite{Raf}.

\subsection{The case with source $\rho_A$}

We had mentioned before, in the context of the Type IIB flow, the comparative lack of development of the theory of geometric flows with singular data compared to the theory for smooth data. The situation is worse for the Type IIA flow, because the source $\rho_A$ is required to satisfy a condition that depends on the form $\varphi$, namely that it be the Poincare dual of a linear combination of special Lagrangians. The notion of Lagrangian depends only on the given symplectic structure, but the notion of special Lagrangian \cite{HL} depends on the metric, which is $g_\varphi$ and hence depends on $\varphi$.

\smallskip
Thus, if we approach the Type IIA equation with brane by flow methods, we have to let both $\rho$ and the $3$-form $\varphi$ evolve. This would already be a challenging and interesting problem in itself. Alternatively, when $\rho_A$ has Dirac-type singularities, we can hope for a scenario where $\rho_A$ arises dynamically as the singularities in the case where the flow with no source does not have any smooth stationary point. From this point of view, the problem becomes reminiscent of a free boundary problem. It would certainly be a very exciting development for both theories if methods from the theory of free boundary problems can be brought successfully to bear on the PDE's from string theory.

\subsection{Related flows in symplectic geometry}

\smallskip

For most of the geometric flows in geometric analysis which are weakly parabolic but not strictly parabolic, e.g. the Ricci flow, the Yang-Mills flow, spinor flows, the Type IIA flows, etc.
we can still establish their short-time existence for arbitrary data. However there are some flows which pose a challenge even in this regard. Some notable examples are the following.

\subsubsection{The Hitchin functional and corresponding gradient flow}

On any compact $6$-d manifold $M$, Hitchin introduced the functional on $3$-forms
$$
H(\varphi)={1\over 2}\int_M \varphi\wedge (J_\varphi\varphi)
$$
and showed that $\delta H=0$ is equivalent to $d(J_\varphi\varphi)=0$, and to $J_\varphi$ being integrable. If $M$ is equipped with a symplectic form and $\varphi$ is primitive, then we get a compatible metric, and we can consider the gradient flow of $H(\varphi)$. It turns out that it is given by
$$
\partial_t\varphi=dd^\dagger\varphi
$$
Recall that Bryant's Laplacian flow on $7$-manifolds is given exactly by this formula, so the Hitchin gradient flow can be viewed as the $6$-manifold version of the $7$-manifold Bryant's Laplacian flow. Remarkably, the Hitchin gradient flow can be cast in turn as a limiting version of the Type IIA flow
$$
\partial_t\varphi=d\Lambda d(\star\varphi)
$$

\subsubsection{The dual Ricci flow}

This flow was introduced by T. Fei and the author \cite{FP1}, as the symplectic dual of the Ricci flow on a complex manifold. It is a flow of positive and primitive $3$-forms on a $6$-dimensional symplectic manifold $M$, which can be worked out to be given by
$$
\partial_t\varphi
=d\Lambda d(\log|\varphi|\star\varphi)
$$
which is similar to the Hitchin gradient flow, but this time with a $\log|\varphi|$ inserted.

\subsubsection{The regularized gradient Hitchin flow}

On the other hand, it is not hard to see that, for any $\epsilon>0$, the modified Type IIA flow defined by
$$
\partial_t\varphi
=d\Lambda d(|\varphi|^\epsilon\star\varphi)
$$
does admit short-time existence, and it is tempting that, from the possibly renormalized solutions $\varphi_\epsilon$ to these flows, we can extract a finite regularized limit, which can serve as solution to either the gradient Hitchin flow or the dual Ricci flow.

\smallskip

Analytically, the situation bears some similarity to the Uhlenbeck-Yau proof \cite{UY} of the Donaldson-Uhlenbeck-Yau theorem on the equivalence between the existence of Hermitian-Einstein metrics and Mumford stability. There they prove that for small $\epsilon>0$, the regularized equation
$$
\Lambda F_\epsilon-\mu I=-\epsilon \log h_\epsilon
$$
always admits a solution. Here $H_\epsilon$ is the unknown metric on a holomorphic vector bundle $E\to (X,g_{\bar kj})$, $F=-\bar\partial(H_\epsilon^{-1}\partial H_\epsilon)$, and $h=H_0^{-1}H_\epsilon$ is the corresponding endomorphism, with $H_0$ a fixed reference metric. The issue then is whether $h_\epsilon$ is uniformly bounded as $\epsilon\to 0$, in which case one can let $\epsilon$ go to $0$ and obtain a solution of the Hermitian-Einstein equation
$$
\Lambda F-\mu I=0
$$
or else, from any subsequence $h_\epsilon$ with ${\rm Tr} h_\epsilon\to\infty$ construct a destabilizing sheaf violating the Mumford stability condition.

\medskip
In this light, it may be the case that a given initial condition may have to satisfy some stability condition for the flow to admit even short-time existence, and the underlying manifold may have to satisfy a stability condition for the short-time existence for any initial data. Of course, these considerations are at the moment quite preliminary, and there is no yet no other strongly supporting evidence for this scenario. Nevertheless, these above flows are geometrically very appealing, and this scenario would seem to be entirely new. So the uncertainty around the approach may be well worth the possible pay-off.

\section{Flows from the Heterotic string}
\setcounter{equation}{0}

Shortly after the identification by Candelas, Horowitz, Strominger, and Witten \cite{CHSW}
of Calabi-Yau $3$-folds as supersymmetric vacua for the heterotic string, Hull \cite{Hull2} and Strominger \cite{S} independently proposed a more general, still supersymmetric, solution allowing metrics with torsion, as long as a key anomaly cancellation is still implemented. Their equations can be described as follows.

\subsection{The Hull-Strominger system}

Let $M$ be
 a compact complex 3-manifold, equipped with a nowhere vanishing holomorphic $(3,0)$-form $\Omega$, and let $E\to M$ be a holomorphic vector bundle with $c_1(E)=0$. Then the Hull-Strominger system is the following system of equations for a Hermitian metric $\o=ig_{\bar kj}dz^j\wedge d\bar z^k$ on $T^{1,0}(M)$ and a Hermitian metric $H_{\bar\alpha\beta}$ on $E$,
 \bea
&& i\p\bar\p\o-{\alpha'\over 4}{\rm Tr}(Rm\wedge Rm-F\wedge F)=0,
\quad
d(\|\Omega\|_\o\o^2)=0\\
&&
F^{2,0}=F^{0,2}=0,
\quad  g^{j\bar k}F_{\bar kj}=0.
\eea
Here $Rm$ and $F$ are the curvatures of $\omega$ and $H$, viewed as $2$-forms valued in $End(T^{1,0}(M))$ and $End(E)$ respectively. The second equation in the first line above was originally written as an equation on the torsion of the metric $\o$. Its reformulation as $d(|\Omega|_\o\o^2)=0$ is an important insight of Li and Yau \cite{LY}, which makes apparent its cohomological meaning. As mentioned earlier, in the mathematical literature, metrics $\o$ on an $m$-dimensional manifold satisfying $d\o^{m-1}=0$ are said to be {\it balanced} \cite{Michelsohn}, and we shall say that a metric $\o$ is conformally balanced if it is balanced up to a conformal transformation. Note a distinguishing features of the equations for the heterotic string, compared to the other string theories, namely the presence of a gauge field $F$ and the Green-Schwarz anomaly cancellation requirement as expressed in the first line of the above system. When $\alpha'=0$, this equation reduces to the Type IIB equation with no source.

\smallskip
Calabi-Yau manifolds can be recognized as special solutions of the Hull-Strominger system with $E=T^{1,0}(M)$, $\o=H$ K\"ahler, and $Ric(\o)=0$. Solutions with singularities were given in \cite{S}, but it was Fu and Yau \cite{FY1, FY2} who found the first smooth non-perturbative, non-K\"ahler solutions, given by toric fibrations over $K2$ surfaces, and who brought to light the deep geometric structure of Hull-Strominger systems. Many other solutions have been found since \cite{FGV, GF2, CPY1}. More generally, it has been suggested by Yau that the Hull-Strominger systems should be viewed as defining canonical metrics for non-K\"ahler complex manifolds with $c_1(M)=0$, and the first steps applying them towards Reid's fantasy have very recently been carried out in \cite{CPY, CGPY}. 

\subsection{Anomaly flows}
As in the previous cases with the equations of the Type IIA and Type IIB strings, a major difficulty in solving the Hull-Strominger system in general is to implement the conformally balanced condition constraint. The same strategy of using a flow that preserves this constraint was advocated in \cite{PPZ1,PPZ2}, and the following flow of pairs $(\omega,H_{\bar\alpha\beta})$ was introduced,
$$
\partial_t(|\Omega|_\omega\omega^2)=
i\partial\bar\partial \omega-{\alpha'\over 4}{\rm Tr}(Rm\wedge Rm-F\wedge F),
\quad
H^{-1}\partial_tH={\omega^2\wedge F\over\omega^3}
\eqno(8)
$$
with initial data $(\omega_0,H_0)$ satisfying the conformally balanced condition $d(|\Omega|_{\omega_0}\omega_0^2)=0$. The flow was called the Anomaly flow, as its stationary points would satisfy the Green-Schwarz Anomaly cancellation equation.
For fixed $\omega$, the flow in $H$ is the well-known Donaldson heat flow, so the novelty in the Anomaly flow lies in the flow for $\omega$, and its coupling to the Donaldson heat flow.
The short-time existence for the Anomaly flow can be established \cite{PPZ1} using Hamilton's Nash-Moser method for weakly parabolic systems \cite{Ha82}.
Using the Anomaly flow, we can recapture [87] the solution of the Hull-Strominger system found by Fu and Yau \cite{FY1, FY2}. The indications are that it will be a powerful tool in the future, and its behavior on unimodular groups and nilmanifolds have been worked out in \cite{PPZ5} and \cite{PU} respectively.

\subsection{The Hull-Strominger-Ivanov system}

Initially, it was believed that the system of equations proposed by Hull and Strominger would automatically imply the equations of motion, but it was subsequently shown by Ivanov and others \cite{Ivanov, FIUV, FIUV2} that this was not the case. Rather, for the field equations to be satisfied, the curvature $Rm$ in the Anomaly cancellation equation must be the curvature $R_\na$ of an $SU(3)$ holomorphic instanton $\na$ on the tangent bundle of $M$. It is not required to be the curvature of the Chern unitary connection defined by $\o$.
We shall refer to the Hull-Strominger system, with the modifications suggested by Ivanov et al, as the Hull-Strominger-Ivanov, or HSI, system. A formulation of this system, which may lead to existence and stability conditions for the existence of solutions, has been proposed by Garcia-Fernandez and Molina \cite{GFMa}, and a Futaki invariant introduced in \cite{GFM}. A formulation in terms of generalized geometry has been proposed in \cite{GFRST}.

\smallskip
Analytically, the HSI system seems also considerably simpler than the original Hull-Strominger system. Indeed, once we fix the holomorphic structure of the $SU(3)$ instanton, its curvature $R_\na$ can be viewed as known,  and there is no more in the Anomaly cancellation equation a quadratic term in the Riemannian curvature $Rm$ of $\o$.
It is straightforward to adapt the geometric flows introduced earlier for the Hull-Strominger system to the HSI system, resulting in an even simpler flow.

\section{Flows from $11$-dimensional supergravity}
\setcounter{equation}{0}

As mentioned above, it is believed that string theories can be themselves unified into M Theory, a theory which has not been fully constructed as of today, and which is mainly understood through its various limits \cite{BBS, HW, To, W}. One of these limits is $11$-dimensional supergravity \cite{CJS}, whose tensor fields are actually very simple: they consist of a Lorentz metric $G_{MN}$ and a closed $4$-form $F=dA$, and the action with all spinor fields set to $0$ is given by 
$$
I(G_{MN},F)
=
\int d^{11}x \sqrt{-G}
(R-{1\over 2}|F|^2)
-{1\over 6}\int A\wedge F\wedge F
$$
The field equations are given by the critical points of the action functional, and given explicitly by
$$
d\star F={1\over 2}F\wedge F,
\qquad
R_{MN}={1\over 12}F_{MPQR}F_N{}^{PQR}.
$$
The supersymmetric solutions are the solutions which admit a non-vanishing spinor $\xi$ satisfying
\bea
D_M\xi\equiv
\na_M\xi-{1\over 288}F_{ABCD}(\Gamma^{ABCD}{}_M+8\,\Gamma^{ABC}\delta^D{}_M)\xi
=0,
\eea
i.e. spinors which are covariantly constant with respect to the connection $D_M$, obtained by twisting the Levi-Civita connection by the flux $F$.

\subsection{Early solutions}

Some early solutions were found with the Ansatz $M^{11}=M^4\times M^7$, where $M^4$ is a Lorentz 4-manifold and $M^7$ a Riemannian manifold with metrics $g_4$ and $g_7$ respectively. Setting $F=c Vol_4$ where $Vol_4$ is the volume form on $M^4$ reduces the field equations to
$$
(Ric_4)_{ij}=-{c^2\over 3}(g_4)_{ij},\qquad
(Ric_7)_{ij}={c^2\over 6}(g_7)_{ij}$$
i.e., $M^4$ and $M^7$ are Einstein manifolds with negative and positive scalar curvatures respectively. These are the Freund-Rubin solutions \cite{FR}, which include $AdS^4\times S^7$.
More sophisticated solutions can be found with other ansatz for $F$, e.g. $F=cVol_4+\psi$ for suitable $\psi$, leading to nearly $G_2$ manifolds (see e.g. \cite{E, DH, PvN, PW} and many others).

\medskip

For us solutions of particular interest are obtained by setting $M^{11}=M^3\times M^8$, and
{$$
g_{11}=e^{2A}g_3+g_8,
\qquad
F=Vol_3\wedge df$$}
where $g_3$ is a Lorentz metric on $M^3$, $g_8$ is a Riemannian metric on $M^8$, and $(A,f)$ are smooth functions on $M_8$. The now well-known solution of Duff-Stelle \cite{Du, DS} is then obtained by assuming the flatness of $g_3$, the conformal flatness of $g_8$, the radial dependence of $A,f$, and supersymmetry.
In \cite{FGP1}, a characterization is provided for when the data $(g_3,g_8,A,f)$ gives rise to a supersymmetric solution:

\begin{theorem}
{\rm \cite{FGP1}}
The data $(g_3,g_8,A,f)$ is a supersymmetric solution to 11-dimensional supergravity equation if and only if

{\rm(a)} $g_3$ is flat;

{\rm (b)} $\bar g_8 := e^Ag_8$ is a Ricci-flat metric admitting covariantly constant spinors with respect to the Levi-Civita connection;

{\rm (c)} $e^{-3A}$ is a harmonic function on $(M_8,g_8)$ with respect to the metric $\bar g_8$; 

{\rm (d}) $df = \pm d(e^{3A})$.
\end{theorem}

Applying the classic results of S.Y. Cheng and P. Li \cite{CL} on lower bounds for Green's functions together with new complete K\"ahler-Ricci flat metrics on ${\bf C}^4$ constructed by Szekelyhidi \cite{Sze}, Conlon-Rochon \cite{CR}, and Li \cite{Li}, we obtain in this manner complete supersymmetric multi-membrane solutions, including many not known before in the physics literature.

\smallskip
If we allow solutions of $11d$ supergravity which are not supersymmetric, then there are many more solutions. For example, it is shown in \cite{FGP1} that the Duff-Stelle solution can be naturally be imbedded in a 5-dimensional family of solutions, in which it is the only supersymmetric solution.

\subsection{Parabolic dimensional reductions of $11d$ supergravity}

More general solutions of $11$-dimensional supergravity will ultimately have to be found by solving partial differential equations. As a starting point, we shall consider the following flow of pairs $(G_{MN}, F)$,
\bea
\label{11d}
\p_tF&=&-\Box_GF-{\sigma\over 2}d\star (F\wedge F)\nonumber\\
\p_tG_{MN}&=&-2R_{MN}+F_{MN}^2-{1\over 3}|F|^2G_{MN}
\eea
where $\Box_G$ is the Hodge-Laplacian on $4$-forms with respect to $G_{MN}$, $F_{MN}^2={1\over 3!}F_{MABC}F_N{}^{ABC}$, and we have introduced a parameter $\sigma=\pm1$, depending on whether the metric $G_{MN}$ is respectively Lorentz or Euclidian. Clearly the solutions of the field equations of $11d$ supergravity must be stationary points of this flow. Conversely, the stationary points of the flow can readily seen to satisfy the field equations if $M^{11}$ is Euclidian, compact, and $H^3(M^{11},{\bf R})=0$. More generally, this remains to be the case, if $M^{11}$ is assumed to have no non-trivial 
$3$-forms which are both closed and co-closed.

\medskip
We shall consider two situations when the above system can be reduced to a parabolic system. The first, immediate situation, arises by considering manifolds $M^{11}$ with Euclidian signature, and $\sigma=-1$. Even though $11$-dimensional supergravity requires a Lorentz signature, it may be expected that this Riemannian version would prove to be of mathematical interest, just as in the case of the Euclidian Yang-Mills theory. In this case, the following theorem was established in \cite{FGP2}:

\begin{theorem}
{\rm \cite{FGP2}}
Consider the flow (\ref{11d}) of a pair $(G_{MN},F)$ with $G_{MN}$ a Euclidian metric, $F$ a closed $4$-form on an $11d$ Riemannian manifold $M^{11}$, and $\sigma=-1$. Then we have

{\rm (a)} The form $F$ remains closed as long as the flow exists;

{\rm (b)} The flow exists for at least a time interval $[0,T_0)$ with $T_0>0$;

{\rm (c)} If $T<\infty$ is the maximum existence time of the flow, then
\bea
{\rm limsup}_{t\to T^-}{\rm sup}_{M^{11}}(|Rm|+|F|)=\infty.
\eea
\end{theorem}

Next, we consider the case of $M^{11}$ of Lorentz signature. In this case, the flow (\ref{11d}) will be hyperbolic. However, we shall identify dimensional reductions of the flow where the Lorentz components are static, and the flow reduces to a parabolic flow.
Thus
we consider space-times and field configurations of the form $M^{11}=M^{1,p}\times M^{10-p}$,
\bea
\label{11dansatz}
G=e^{2A}g_{1,p}+g,
\quad F=dVol_g\wedge\beta+\Psi
\eea
where $\beta$ and $\Psi$ are now respectively a $(3-p)$-form and a $4$-form on $M^{10-p}$, both of whom are assumed to be closed. The metric $g_{1,p}$ is a Lorentz and Einstein metric on $M^{1,p}$,
\bea
Ric(g_{1,p})=\lambda g_{1,p},
\eea
$A$ is a scalar function on $M^{10-p}$, and the metric $g=g_{ij}$ is a Riemannian metric on $M^{10-p}$.
The dimension $p$ is allowed to take any value between 0 and 10, and $\beta=0$ if $p\geq 4$. 

\medskip
For static $g_{1,p}$, the geometry of $M^{11}$ is characterized by the $4$-tuple $(g_{ij}, A,\beta, \Psi)$. We shall consider the following flow of such $4$-tuples:
\bea
\label{11d-reduced}
\p_tg_{ij}&=&-2R_{ij}+2(p+1)(\na_i\na_jA+\na_iA\na_jA)\nonumber\\
&&\quad-e^{2(p+1)A}\beta_{ij}^2+\Psi_{ij}^2
-{1\over 3}(|\Psi|^2-e^{-2(p+1)A}|\beta|^2)g_{ij}\nonumber\\
\p_tA&=&\Delta A+{1\over 4}(p+1)|\na A|^2-{1\over 3}e^{-2(p+1)A}|\beta|^2
-{1\over 6}|\Psi|^2-\lambda e^{-2A}
\nonumber\\
\p_t\beta&=&-\Box\beta+
(-1)^{p+1}(p+1)d\star (dA\wedge \star\beta)
-(-1)^p{p+1\over 2}e^{(p+1)A}dA\wedge \star(\Psi\wedge\Psi)
\nonumber\\
&=&\quad-(-1)^p{1\over 2}e^{(p+1)A}d\star(\Psi\wedge\Psi)
\nonumber\\
\p_t\Psi
&=&-\Box\Psi+(-1)^p(p+1)d\star(dA\wedge \Psi)-(p+1)e^{-(p+1)A}dA\wedge \star(\beta\wedge\Psi)\nonumber\\
&&\quad
e^{-(p+1)A}d\star(\beta\wedge\Psi)
\eea
where the norms, Hodge star operator, and covariant derivatives are all taken with respect to the metric $g_{ij}$ on $M^{10-p}$. The following theorem was established in \cite{FGP2}:

\begin{theorem} 
{\rm\cite{FGP2}}
Let $M^{11}$ be a Lorentz $11$-dimensional manifold, equipped with a metric $G_{MN}$ and a closed $4$-form $F$ given by the expression (\ref{11dansatz}) in terms of the $4$-tuple $(g_{ij},A,\beta,\Psi)$ on a Riemannian manifold $M^{10-p}$.

{\rm (a)} If the $4$-tuple evolves by the above flow (\ref{11d-reduced}), then the pair $(G_{MN},F)$ on $M^{11}$ evolved by the original flow (\ref{11d}) with $\sigma=1$.

{\rm (b)} The forms $\beta,\Psi$ and $F$ remain closed along the flow, and the above ansatz is preserved;

{\rm (c)} The flow exists at least for some time interval $[0,T_0)$ with $T_0>0$. If $T<\infty$ is the maximum time of existence of the flow, then
\bea
{\rm limsup}_{t\to T^-}{\rm sup}_{M^{10-p}}
(|Rm|+|A|+|\beta|+|\Psi|)=\infty
\eea
\end{theorem}

\medskip
Similar reductions should be possible for $5$-dimensional supergravity. By construction, the stationary points of these flows will satisfy the supergravity field equations. In view of the above stated goal of finding spaces satisfying both supersymmetry and the field equations, it is natural to examine these flows on spaces which are at least known to admit a non-vanishing spinor. The same reductions of $G$ structure that are anticipated to be helpful for the spinor flows of $(g_{ab},\psi, H,\varphi)$ described in the next section \S 6 should also be instrumental in the parabolic reductions flows of $(g_{ij},A,\beta, \Psi)$ described in above. It should be very instructive to have a comparison of the two types of flow under this assumption of existence of a non-vanishing spinor field.

\section{Flows of spinor fields}
\setcounter{equation}{0}

While this is not apparent at first sight, 
all the equations for unified string theories which we have written down so far, except for $11d$  supergravity, originate from the requirement of (local) supersymmetry. Local supersymmetry can be viewed as an analogue of reparametrization invariance, whose transformations are generated by a spinor field instead of a vector field. Space-time is required to be invariant under such transformations, which reduces to the condition that the generating spinor field be a non-zero field which is covariantly constant with respect to a suitable connection possibly with flux. Thus the common underlying thread in all the equations which we have written down so far, except for $11d$ supergravity, is the existence of such a covariantly constant spinor field, and the equations themselves reflect the consequences of the existence of this spinor field on the tensor (``bosonic") fields of the theory, while not involving directly the spinor (``fermionic") field themselves.

\medskip
This suggests that an alternative, and complementary, approach to the equations of unified string theories would be to examine conditions under which such a covariantly constant spinor field can exist, as they are necessary for supersymmetry. Such an approach, especially if it was based on geometric flows, may have its own advantages: already we have observed that, in the case of $11d$ supergravity, this is a natural approach. But perhaps more important, 
it is a fundamental aspect of unified string theories, which we had not yet addressed except in special cases in \S 3, that they are mutually related by dualities. Such dualities are highly non-trivial, and usually hard to detect in theories which have each their own formulation. It can be hoped that this task would become easier if we restrict to the space of supersymmetric theories, and more specifically to the moduli space of theories admitting a covariantly constant spinor. Such a moduli space and Kuranishi deformation theory, especially with fluxes, may well be of considerable geometric interest in its own right.

\smallskip
We observe that the existence of a covariantly constant spinor is well known in mathematics to be characteristic of reduced holonomy and special geometry (see e.g. Berger \cite{Be}, Lichnerowicz \cite{Li}, Joyce \cite{J} and others).
An approach using geometric flows had been initiated by Ammann, Weiss, and Witt \cite{AWW}, based on the gradient flow of the energy functional $I(g_{ij},\psi)=\|\na\psi\|_{L^2}^2$, about which we shall say more below. What the physics of unified string theories has done is to infuse a new interest in such questions, with the key additional complication of the possible presence of flux.

\subsection{Spinors, connections, and supersymmetry}

We begin by recalling the basic set-up for spinor fields. Let $M$ be a compact spin manifold of dimension $n$, and $S\to M$ the bundle of spinors on $M$. We always work in local trivializations of the spin bundle $S$, which project on local trivializations of the bundle of orthonormal frames, given by a choice of frame $\{e_a\}
=\{e_a{}^j\}$ in a local coordinate system $\{x^j\}$, $1\leq j\leq n$. The space of spinors is given then pointwise by a vector space $V$, so that spinor fields are just locally functions $\psi$ valued in $V$. A defining property of the vector space $V$ is that it carries a representation of the Clifford algebra, namely the algebra generated by Dirac matrices $\{\gamma^a\}$. We shall adopt the following convention for Dirac matrices,
\bea
\gamma^a\gamma^b+\gamma^b\gamma^a=2\,\delta^{ab}.
\eea
A basis for the endomorphisms of $V$ is given by the anti-symmetrized products
\bea
\gamma^{b_1\cdots b_p}={1\over p!}\sum_{\sigma}(-1)^{|\sigma|} \gamma^{\sigma_1}\cdots\gamma^{\sigma_p}
\eea
where $\sigma:(1,\cdots,p)\to(\sigma_1,\cdots,\sigma_p)$ runs over all permutations of $(1,\cdots,p)$, and $|\sigma|$ is its signature. The spin bundle $S$ admits a natural connection $\na$, namely the spin connection given by
\bea
\na_j\psi=\p_j\psi-{1\over 4}\o_{jab}\gamma^{ab}\psi
\eea
where $\o_{jab}$ is the Levi-Civita connection.

\smallskip

A key feature of all the limits of M Theory is that they automatically incorporate a supergravity theory, that is, a theory which incorporates gravitation and is supersymmetric.
In particular, their spectrum includes the supersymmetric partner of the metric $g_{ij}$, which is a spinor-valued form $\chi_i=\{\chi_i{}^\alpha\}\in T_*(M)\otimes S$ called the gravitino field.
As mentioned above, the generator of supersymmetry transformations is a spinor field $\psi$, and in analogy with generators of diffeomorphisms which are vector fields $V^j$ and act on metrics by $\delta g_{ij}=\na_iV_j+\na_jV_i$, supersymmetry transformations act on the gravitino by
\bea
\delta \chi_i=\na_i\psi+\cdots.
\eea
Here the right hand side is a connection on $\psi$, just as the right hand side in diffeomorphisms is a connection on $V^j$, and the dots $\cdots$ indicate contributions to this connection from the other fields in the theory. We shall generically refer to them as ``flux" terms. Since the endomorphisms of $V$ are generated by the Clifford algebra, it follows that the flux terms must be given by a polynomial in the Dirac matrices $\gamma^a$, whose coefficients depend on the other fields of the theory,
\bea
\sum_p\sum_{b_1,\cdots,b_p}H_{Mb_1\cdots b_p}^{(p)}\gamma^{b_1\cdots{b_p}}\psi
\eea
The cases of particular importance for $M$ theory are when the field $H_{MN_1\cdots N_p}^{(p)}$ is a $(p+1)$-form. The case of a $3$-form is responsible for the torsion $H$ in the equations for the Type I, Type II, and heterotic string, and the case of a $4$-form results in the field $F_4$  in $11$-dimensional supergravity, both discussed earlier. 

\smallskip
Thus for the supergravity theories which we consider and for our purposes, the flux terms 
can be taken to be in a very specific form. Fix a positive integer $p$, and let $H$ be a $p$-form. Then the covariant derivatives on spinors which we consider are of the form
\bea
\nabla_a^H\psi=\nabla_a\psi+(\lambda|H)_a
\eea
where for each vector $\lambda=(\lambda_1,\lambda_2)\in {\bf R}^2$, we have set
\bea
(\lambda|H)_a\psi=\lambda_1H_{ab_1\cdots b_{k-1}}\gamma^{b_1\cdots b_p}\psi+\lambda_2H_{b_1\cdots b_k}\gamma^{ab_1\cdots b_k}\psi,
\eea
which is a $1$-form valued in the space $End(S)$ of endomorphisms of spinors. Note that the two terms in $(\lambda|H)_a$ are either Hermitian or anti-Hermitian, depending on the parity of $p$, but not both. 

\smallskip
We set then $\na_a^H\psi=\na_a\psi+(\lambda|H)_a\psi$, and we are particularly interested in determining when a compact spin manifold $M$ admits a nowhere vanishing spinor field $\psi$ with $\na^H\psi=0$.
As mentioned earlier, in the special case where $H$ is assumed to be $0$, one can consider the gradient flow of the functional
\bea
I(g_{ij},\psi)=\int_X |\na\psi|^2\sqrt g dx
\eea
subject to the constraint $|\psi|^2=1$. This flow was introduced by Ammann, Weiss, and Witt \cite{AWW}, who established its short-time existence. Shi-type estimates for it were subsequently established by He and Wang \cite{HW}.
However, the possibility of a non-vanishing flux $H$ leads to considerable new complications, not least due to the fact that connections with fluxes are usually not unitary, and that there is no obvious normalization of $\psi$ that would prevent it from sliding to $0$. 

\smallskip
In \cite{CP}, T. Collins and the author proposed to address this problem by introducing a dynamic normalization condition
\bea
\label{normalization-varphi}
e^{-2\varphi}|\psi|^2=1
\eea
and considering flows of the $4$-tuple $(g_{ij},\psi,H,\varphi)$. While the arguments apply with some generality, they do depend on relations between the dimension $n$ of the spin manifold $M$ and the rank $p$ of the flux $H$. In the following we shall assume that $3p=n+1$, which correspond to the cases of $5d$ and $11d$ supergravity. With this assumption, we consider the following flow of the $4$-tuple $(g_{ij},\psi,H,\varphi)$,
 \bea
 \label{spinorflow}
\partial_tg_{ab}&=&{1\over 4}e^{-2\varphi}{\rm Re}\,Q_{\{ab\}}^H+
 {1\over 2}\<\nabla_{\{a}^H\psi,\nabla_{b\}}^H\psi\>-{1\over 4}g_{ab}|\nabla^H\psi|^2)
 \nonumber\\
\partial_tH&=&-(d^\dagger dH+d(d^\dagger H-c\star (H\wedge H))-L(H))
\nonumber\\
\partial_t\varphi&=&\sum_ae^{-2\varphi}
(\na_a-2e^{-2\varphi}\<h_a\psi,\psi\>)
\{e^{2\varphi}(\na_a\varphi+e^{-2\varphi}(\na_a\varphi+e^{-2\varphi}
\<h_a\psi,\psi\>)\}\nonumber\\
\partial_t\psi&=&-\sum_a(\nabla_a^H)^\dagger\nabla_a^H\psi+c(t)\psi
 \eea
The various quantities entering the formula are defined as follows: 
brackets $\{ab\}$ denote symmetrization with respect to $a$ and $b$; $h_a={1\over 2}
((\lambda|H)_a+(\lambda|H)_a^\dagger)$ is the Hermitian part of the operator 
$(\lambda|H)_a$; $c$ is a given constant; 
$L(H)$ is the $p$-form defined by
$$
\<L(H),\beta\>={1\over 2}\<c\star (\beta\wedge H+H\wedge\beta),d^\dagger H-c\star (H\wedge H)\>,
$$
the coefficients $c(t)$ and tensor $Q_{ab}$ are defined by
\bea
&c(t)
=e^{-2\varphi}|\nabla^H\psi|^2-2e^{-2\varphi}
\sum_a \<h_a\psi,\psi\>(
(\na_a\varphi+e^{-2\varphi}\<h_a\psi,\psi\>)\nonumber\\
&
Q_{ab}^H=\nabla_c\<\gamma^{ca}\psi,\nabla_b^H\psi\>.
\nonumber
\eea

\begin{theorem}
{\rm \cite{CP}}
Consider the flow (\ref{spinorflow}) of $4$-tuples $(g_{ij},\psi, H,\varphi)$ on a spin manifold $M$ of dimension $n$, where $H$ is a $p$-form, and $3p=n+1$. Then it preserves the normalization condition (\ref{normalization-varphi}). Furthermore,

{\rm (a)} At any stationary point, the spinor field is covariantly constant with respect to the connection $\na^H$, and the flux $H$ satisfies the equation
\bea
\label{H-equation}
dH=0,
\quad
d^\dagger H=c\star H\wedge H.
\eea

{\rm (b)} The flow (\ref{spinorflow}) is weakly parabolic, and it always exists for some positive time interval, for any smooth initial data satisfying the normalization constraint (\ref{normalization-varphi}).

{\rm (c)} The following Shi-type estimates hold: assume that on the time interval $(0,{\tau\over A})$, we have the following estimates
\bea
|\na\psi|^2, |\na^2\psi|,
|H|^2,|\na H|,|Rm|\leq A,
\quad
\|\varphi\|_{L^\infty}\leq \hat A
\eea
for some constant $\hat A$. Then for any integer $q\geq 0$, there is a constant $C_q$ depending only on $q,n,\hat A$, and an upper bound for $\tau$ so that
\bea
|\na^qRm|+|\na^{q+2}\psi|+|\na^{q+2}e^{2\varphi}|+|\na^{q+1}H|\leq C_q{A\over t^{q/2}}
\eea
uniformly on $[0,{\tau\over A})$.
\end{theorem}

We note that the above flows are suggested by the requirement that, at the stationary points, the flux $H$ satisfies the equations arising from $5d$ and $11d$ supergravity theories. A key new observation is that the requirement that $\psi$ be covariantly constant implies that $\varphi$ must satisfy a specific equation of its own. The choice of the flow for $\varphi$ in (\ref{spinorflow}) is designed so that this equation is satisfied at the stationary points, and so that the whole system is weakly parabolic.

\smallskip
It may also be instructive to indicate here an identity between the form $Q_{\{p,\ell\}}^H$ and the Ricci curvature $R_{p\ell}$, generalizing such an identity due to Ammann, Weiss, Witt \cite{AWW} and He and Wang \cite{HW},
\bea
\Re (Q_{\{p,\ell\}}^H)&=&{1\over 2}e^{2\varphi}R_{p\ell}+
\Re(\<\psi,\gamma^\ell\na_p\slash D\psi\>)-{1\over 2}{\rm Hess}(e^{2\varphi})_{p\ell}
\nonumber\\
&&
-\Re(\<\slash D\psi,\gamma^\ell\na_p\psi\>)+\na_a\Re(\<\gamma^{a\ell}\psi,(\lambda|H)_p\psi\>)
+2\Re\<\na_p\psi,\na_\ell\psi\>
\nonumber
\eea
which relates the flow of the metric $g_{p\ell}$ to the Ricci flow.

\smallskip
The following theorem is more specifically dependent on viewing the spinor flow (\ref{spinorflow}) as a coupled system:

\begin{theorem}
{\rm\cite{CP}}
Assume that we have the following estimates
\bea
|\na\psi|^2,
|\na^2\psi|,|H|^2,|\na H|,|\na^2H|^{2\over 3}\leq A,
\qquad
\|\varphi\|_{L^\infty}\leq \hat A
\eea
on some time interval $[0,\tau)$. Then there exists a constant $C$, depending only on $A,\hat A,n$, an upper bound for $\tau$, and
\bea
\kappa:=\int_M |Rm(0)|^q\sqrt{g(0)} +{\rm Vol}(M,g(0))
\eea
for some $q>{n\over 2}$, so that, uniformly on $[0,\tau)$, we have
\bea
|Rm(t)|_{g(t)}\leq C.
\eea
\end{theorem}

Finally, we would like to mention a different interesting direction, namely parallel spinors on Lorentz manifolds. In particular, it has been shown by Murcia and Shahbazi \cite{MS1, MS2} that the existence of a global parallel spinor on a hyperbolic Lorentz 4-manifold induces a solution of a hyperbolic flow called the ``parallel spinor flow" on a Cauchy hypersurface, and vice versa.

\section{Some further developments}
\setcounter{equation}{0}

We had restricted our discussion to geometric flows motivated by string theories, distinctive features of which have been some new underlying special geometry. But there has been considerable progress as well in recent years on flows with more established underlying notions of geometry. One case which has attracted increasing attention has been $G_2$ flows or flows of balanced metrics. The short-time existence and Shi-type estimates for Bryant's $G_2$ Laplacian flow have been established 
in \cite{BryXu, BV, Ka, LoWe, FFR}, solitons introduced in \cite{FiRaf}, and the convergence in particular cases established in \cite{PS}, building partly on techniques newly introduced in the study of Anomaly flows. Another case is the K\"ahler-Ricci flow \cite{Cao}, especially on questions motivated by the Analytic Minimal Model Program \cite{SoT}. New PDE techniques for $L^\infty$ estimates for complex Monge-Amp\`ere equations \cite{GPT, GuoP} have allowed new estimates on Green's functions no longer dependent on direct conditions on the Ricci curvature \cite{GPS}. These in turn have helped answer some long-standing questions on diameter bounds and the Gromov-Hausdorff convergence of the K\"ahler-Ricci flow, in both cases of  finite-time and long-time solutions \cite{GPSS}.

\smallskip
Looking forward, the fact that practically all the flows from unified string theories can be viewed as the Ricci flow coupled to other fields, suggests that the development of methods to study such flows should be very valuable. There appears to have been relatively few efforts in this direction, except for couplings to scalar fields and harmonic maps (see e.g. \cite{Mue, List, Johne, GHP} and references therein). On the other hand, a coupled system which has been of great importance in geometry and physics is the Hitchin system \cite{Hita} on Riemann surfaces. There has been considerable progress recently on the singularities of the Vafa-Witten \cite{VW} and the Kapustin-Witten \cite{KW} equations, which can be viewed as natural generalizations of Hitchin systems to higher dimensions. See \cite{Taubes, Ta}. 
Similar results for coupled Ricci flows could be truly instrumental for the flows considered here.

\bigskip
\noindent
{\bf Acknowledgements} The author is very grateful to his collaborators, Tristan Collins, Teng Fei, Bin Guo, Sebastien Picard, and Xiangwen Zhang, for joint efforts on the works presented here. He would like especially to thank Li-Sheng Tseng for describing to him his joint works with Shing-Tung Yau on supersymmetry and cohomology, which influenced him greatly. Parts of the material in this survey have also been presented in lectures at Harvard University, the University of California at Irvine, the University of Torino, the Laboratoire de Mathematiques at Orsay, Hamburg University, the University of Connecticut at Storrs, and at the conference on ``Complex Analysis and Dynamics" at Portoroz, Slovenia.

\bigskip

email address: phong@math.columbia.edu

Department of Mathematics

Columbia University, New York, NY 10027


\begin{thebibliography}{99}

{\footnotesize

\bibitem{AWW} B. Ammann, H. Weiss, and F. Witt,
{\it A spinorial energy functional: critical points and gradient flow}, Math. Ann. 365 (2016) 1559-1602.

\bibitem{BBS} K. Becker, M. Becker, and J.H. Schwarz,
{\it String theory and M-Theory: A modern Introduction},
Cambridge University Press, 2007.

\bibitem{BV} L. Bedulli and L. Vezzoni, {\it A parabolic flow of balanced metrics}. J. Reine Angew. Math. 723 (2017), 79-99.

\bibitem{BVa} L. Bedulli and L. Vezzoni, {\it Stability of flows of closed forms}, Adv. Math. 364 (2020) 107030, 29 pp.

\bibitem{BVb} L. Bedulli and L. Vezzoni,
{\it On the stability of the anomaly flow}, Math. Res. Lett. 29 (2022) No. 2, 323-338.


\bibitem{BI} D.E. Blair and S. Ianus, ``Critical associated metrics on symplectic manifolds", Nonlinear Problems in Geometry 23-29, Contemp. Math. 51, AMS, 1986.

\bibitem{Bry} R. Bryant, ``Some remarks on $G_2$ stuctures", Proc. of G\"okova Geometry-Topology Conference 2005, 75-109, International Press, 2006.

\bibitem{BryXu} R. Bryant and F. Xu, {\it Laplacian flow for closed G2-structures: short time behavior}, arXiv:1101.2004

\bibitem{Be} M. Berger,
{\it Sur les groupes d’holonomie homog\`enes de vari\'et\'es à connexion affine et des variétés riemanniennes}
Bulletin de la S. M. F., tome 83 (1955) 279-330.

\bibitem{CHSW} P. Candelas, G. Horowitz, A. Strominger and E. Witten, {\it Vacuum configurations for superstrings}, Nucl. Phys. B 258 (1985) 46-74.

\bibitem{Cao} H.D. Cao, {\it Deformation of K\"ahler metrics to K\"ahler-Einstein metrics on compact K\"ahler manifolds}, Invent. Math. 81 (1985), no. 2, 359-372.

\bibitem{CaK} H.D. Cao and J. Keller,
{\it On the Calabi problem: a finite-dimensional approach}, J. Eur. Math. Soc. (JEMS) 15 (2013), no. 3, 1033-1065.

\bibitem{Caoetal} J. Cao, P. Graf, P. Naumann, M. Paun, T. Peternell, and X. Wu,
{\it Hermite-Einstein metrics in singular settings}, arXiv: 2303.08773.

\bibitem{Cardoso} G.L. Cardoso, G. Curio, G. Dall'Agara, D. L\"ust, P. Manousselis, and G. Zoupanos,
{\it Non-K\"ahler string backgrounds and their five torsion classes}, hep-th/0211118.

\bibitem{CL} 
S.Y. Cheng and P. Li, {\it Heat kernel estimates and lower bounds of eigenvalues}, Comm.
Math. Helvetici 56 (1981), 327-338.

\bibitem{CS} S. Chiossi and S. Salamon, {\it The intrinsic torsion of $SU(3)$ and $G_2$ structures}, math.DG/0202282, in {\it Differential Geometry}, Valencia, 2001, pp. 115-133.

\bibitem{CHT} T. Collins, T. Hisamoto, and R. Takahashi, {\it  The inverse Monge-Amp\`ere flow and applications to K\"ahler-Einstein metrics}, J. Differential Geom. 120 (2022) no. 1, 51-95.


\bibitem{CGPY} T. Collins, S. Gukov, S. Picard, and S.T. Yau,
{\it Special Lagrangian cycles and Calabi-Yau transitions}, arXiv: 2111.10355.

\bibitem{CP} T. Collins and D.H. Phong,
{\it Spinor flows with flux, I: short-time existence and smoothing estimates}, arXiv: 2212.00814.

\bibitem{CPY} T. Collins, S. Picard, and S.T. Yau,
{\it Stability of the tangent bundle through conifold transitions}, to appear in Comm. Pure Appl. Math. 

\bibitem{CPY1} T. Collins, S. Picard, and S.T. Yau,
{\it The Strominger system in the square of a K\"ahler class},
arXiv: 2211.03784.



\bibitem{CR} R. Conlon and F. Rochon, {\it New examples of complete Calabi-Yau metrics on $C^n$ for $n\geq 3$},
arXiv: 1705.08788

\bibitem{CJS} E. Cremmer, B. Julia, and J. Scherk, {\it Supergravity theory in 11 dimensions}, Phys. Lett. B
76 (1978), 409-412.

\bibitem{deBT} P. de Bartholomeis and A. Tomassini,
{\it On the Maslov index of Lagrangian submanifolds of generalized Calabi-Yau manifolds}, Internat. J. Math. 17 (2006) no. 8, 921-947.


\bibitem{DH} E. D'Hoker, 
{\it Exact M-Theory Solutions, Integrable Systems, and Superalgebras}, SIGMA 11
(2015) 609-628.

\bibitem{D} S.K. Donaldson, {\it Anti self-dual Yang-Mills connections over complex algebraic surfaces and stable vector bundles}, Proc. London Math. Soc. (3) 50 (1985), no.1, 1-26.

\bibitem{Du} M.J. Duff, 
{\it The world in eleven dimensions: supergravity, supermembranes and M-theory},
(1999) Studies in High Energy Physics and Cosmology, IOP Publishing, Taylor \& Francis.

\bibitem{DS} M.J. Duff and K.S. Stelle, {\it Multi-membrane solutions of D=11 supergravity}, Phys. Lett. B
253 (1991), 113-118.

\bibitem{E} F. Englert,
{\it Spontaneous compactification of eleven-dimensional supergravity}, Phys. Lett. B
119 (1982), 339 -342.

\bibitem{FGP1} T. Fei, B. Guo, and D.H. Phong, {\it A geometric construction of solutions to 11d supergravity}, Comm. Math. Phys. 369 (2019) no. 2, 811-836.

\bibitem{FGP2} T. Fei, B. Guo, and D.H. Phong, {\it Parabolic dimensional reductions of 11d supergravity}, Analysis \& PDE Vol 14 (2021) No. 5, 1333-1361.

\bibitem{FP} T. Fei and D.H. Phong, {\it Unification of the K\"ahler-Ricci and Anomaly flows}, 
Surveys Diff. Geom. XXIII (2020) 89-103.

\bibitem{FP1} T. Fei and D.H. Phong, {\it Symplectic geometric flows}, arXiv: 2111.14048, to appear in PAMQ, issue in honor of V. Guillemin.

\bibitem{FPPZb} T. Fei, D.H. Phong, S. Picard, and X.W. Zhang,
{\it Estimates for geometric flows for the Type IIB string}, Math. Ann. 382 (2022) 
1935-1955.

\bibitem{FPPZ} T. Fei, D.H. Phong, S. Picard, and X.W. Zhang, {\it Geometric flows for the Type IIA string}, Cambridge J. Math. Vol. 9 (2021) No. 3, 693-807.

\bibitem{FPPZa} T. Fei, D.H. Phong, S. Picard, and X.W. Zhang, {\it Stability of the Type IIA flow and its applications in symplectic geometry}, arXiv: 2112.15580, to appear in Comm. Analysis and Geometry.

\bibitem{FFR} M. Fernandez, A. Fino, and A. Raffero,
{\it On G2-structures, special metrics and related flows}
arXiv:1810.07587.

\bibitem{FIUV} M. Fernandez, S. Ivanov, L. Ugarte and R. Villacampa, {\it Non-K\"ahler heterotic string compactifications with non-zero fluxes and constant dilaton}, Comm. Math. Phys. 288 (2009), 677-697.

\bibitem{FIUV2} M. Fernandez, S. Ivanov, L. Ugarte and D. Vassilev, {\it Non-K\"ahler heterotic string solutions with non-zero fluxes and non-constant dilaton}, J. of High Energy Phy., Vol 6, (2014) 1-23.

\bibitem{FGV} A. Fino, G. Grantcharov, and L. Vezzoni,
{\it Solutions to the Hull-Strominger system with torus symmetry}, arXiv: 1901.10322.

\bibitem{FiRaf} A. Fino and A. Raffero,
{\it Remarks on homogeneous solitons of the $G_2$ Laplacian flow},
Comptes Rendus Mathematique 358 (4) (2020) 401-406.

\bibitem{FR} P. Freund and M. Rubin, {\it Dynamics of dimensional reduction}, Phys. Lett. B 97
(1980), 233 – 235.

\bibitem{Fr} T. Friedrich, {\it On types of non-integrable geometries}, math.DG/0205149.

\bibitem{FY1} J.X. Fu and S.T. Yau, {\it The theory of superstring with flux on non-K\"ahler manifolds and the complex Monge-Amp\`ere equation}, J. Differential Geom., Vol 78, Number 3 (2008), 369-428.

\bibitem{FY2} J.X. Fu and S.T. Yau, {\it A Monge-Amp\`ere type equation motivated by string theory}, Comm. Anal. Geom. 15 (2007), no. 1, 29-76.

\bibitem{GF2} M. Garcia-Fernandez,
{\it T-dual solutions of the Hull-Strominger system on non-K\"ahler threefolds},
arXiv:1810.04740.

\bibitem{GFRST} M. Garcia-Fernandez, R. Rubio, C.S. Shahbazi, C. Tipler, {\it Canonical metrics on holomorphic Courant algebroids},
arXiv:1803.01873  [pdf, ps, other]   math.DG hep-th math.AG math.SG

\bibitem{GFM} M. Garcia-Fernandez and R.G. Molina, {\it Futaki invariants and Yau's conjecture on the Hull-Strominger system},
arXiv: 2303.04274.

\bibitem{GFMa} M. Garcia-Fernandez and R.G. Molina,
{\it Harmonic metrics for the Hull-Strominger system and stability},
arXiv: 2301.08236.

\bibitem{Gau} P. Gauduchon, 
{\it Hermitian connections and Dirac operators}, Boll. Un. Mat. Ital. B (7) 11 (1997) no. 2, suppl. 257-288.

\bibitem{GGHPR}
J. Gauntlett, J. Gutowski, C. Hull, S. Pakis, and H. Reall, {\it All supersymmetric solutions of minimal supergravity in five-dimensions}, Classical Quantum Gravity, 20 (2003), 4587. [hep-th/0209114]

\bibitem{GMW}
J. Gauntlett, D. Martelli, and D. Waldram,
{\it Superstrings with Intrinsic Torsion},
arXiv:hep-th/0302158, also doi Phys. Rev. D.69.086002. 


\bibitem{GP}
J. Gauntlett and S. Pakis, {\it The geometry of D = 11 Killing spinors}, arXiv:hep-th/0212008, JHEP 039 (2003), 32 pp.

\bibitem{Gr} M. Grana, R. Minasian, M. Petrini, and A. Tomasiello, {\it Supersymmetric backgrounds
from generalized Calabi-Yau manifolds}, J. High Energy Phys. (2004), no. 08, 046,
arXiv:hep-th/0406137. 

\bibitem{Guanbo} B. Guan, { \it Second order estimates and regularity for fully nonlinear elliptic equations on Riemannian manifolds}, Duke Math. J. 163 (2014), 1491-1524.

\bibitem{GHP} B. Guo, Z. Huang, and D.H. Phong,
{\it Pseudo-locality for a coupled Ricci flow}, Comm. Anal. Geom. 26 (2018) no. 3, 593-626.

\bibitem{GuoP} B. Guo and D.H. Phong,
{\it On $L^\infty$ estimates for fully non-linear partial differential equations},
arXiv: 2204.12549.

\bibitem{GPS} B. Guo, D.H. Phong, and J. Sturm,
{\it Green's functions and complex Monge-Amp\`ere equations}, arXiv: 2204.12549.

\bibitem{GPSS} B. Guo, D.H. Phong, J. Song, and J. Sturm,
{\it Diameter estimates in K\"ahler geometry}, arXiv: 2209.09428.

\bibitem{GPT} B. Guo, D.H. Phong, and F. Tong,
``On $L^\infty$ estimates for complex Monge-Amp\`ere equations", to appear in Ann. of Math. (2023).

\bibitem{Gurrieri} S. Gurrieri, J. Louis, A. Micu, and D. Waldram,
{\it Mirror symmetry in Calabi-Yau compactifications},
arXiv hep-th/0211102.

\bibitem{GMR} J.B. Gutowski, D. Martelli, and H. S. Reall, {\it All Supersymmetric solutions of minimal supergravity in six- dimensions}, Class. Quant. Grav. 20 (2003) 5049–5078, hep-th/0306235

\bibitem{Ha82} R. Hamilton, {\it Three-manifolds with positive Ricci curvature}, J. Differential Geom. 17 (1982) no. 2, 255-306.

\bibitem{Ha95} R. Hamilton, {\it The formation of singularities in the Ricci flow}, Surveys in Differential Geom. Vol. II (1993).

\bibitem{HL} R. Harvey and H.B. Lawson,
{\it calibrated geometries}, Acta Math. 148 (1982) 47-157.

\bibitem{HW} F. He and C. Wang,
{\it Regularity estimates for gradient flow of a spinorial energy functional}, to appear in Math. Res. Lett.

\bibitem{Hita} N.J. Hitchin,
{\it The self-duality equation on a Riemann surface}, Proc. London Math. Soc. 55 (1987) 59-126.

\bibitem{Hitb} N.J. Hitchin,
{\it The geometry of three-forms in six dimensions},
J. Differential Geom. 55 (2000) No. 3, 547-576.


\bibitem{HW} P. Horava and E. Witten, {\it Heterotic and type I string dynamics from eleven dimensions},
Nucl. Phys. 460 (1995) 506.

\bibitem{Hull2}
C. Hull, {\it Compactifications of the Heterotic Superstring}, Phys. Lett. 178 B (1986) 357-364.

\bibitem{Ivanov} S. Ivanov,
{\it Heterotic supersymmetry, anomaly cancellation, and equations of motion},
Phys. Lett. B 685 (2010) 190-196.

\bibitem{Johne} F. Johne, {\it Surgery for extended Ricci flow systems}, Thesis Univ. of T\"ubingen (2019).

\bibitem{J} D. Joyce, {\it Compact Manifolds with Special Holonomy}, Oxford Math. Monogr., Oxford University
Press, Oxford, 2000.

\bibitem{KW} A. Kapustin and E. Witten,
{\it Electric-magnetic duality and the geometric Langlands program},
Comm. Number Theory Physics 1 (2007) 1-236.

\bibitem{Ka} S. Karigiannis, {\it Flows of $G2$-structures, I.} Quarterly J. of Math., 60(4) (2009), 487-522.

\bibitem{K1} N. Klemyatin,
{\it Convergence of Hermitian manifolds and the Type IIB flow},
arXiv:2209.00312.

\bibitem{K2} N. Klemyatin, 
{\it Convergence of Type IIA manifolds and application to the Type IIA flow},
arXiv: 2210.05082.

\bibitem{LW} H.V. L\^e and G.F. Wang, {\it Anti-complexified Ricci flow on compact symplectic manifolds},
J. Reine Angew. Math. 530 (2001) 17-31.

\bibitem{Li} Y. Li, {\it A new complete Calabi-Yau metric on ${\bf C}^3$}, arXiv: 1705.07026.

\bibitem{LY} J. Li and S.T. Yau, {\it The existence of supersymmetric string theory with torsion}, J. Differ. Geom. 70 (2005), no.1, 143-181.

\bibitem{Li} A. Lichnerowicz, 
{\it Théorie globale des connexions et des groupes d'holonomie}, 
Edizioni Cremonese, Roma (1955) 15+282 pp.

\bibitem{List} B. List, {\it Evolution of an extended Ricci flow system},
Comm. Anal. Geom. 16 (2008) no. 5, 1007-1048.

\bibitem{LoWe} J. Lotay and Y. Wei, {\it Laplacian flow for closed G 2 structures: Shi-type estimates, uniqueness and compactness}, Geometric and Functional Analysis, 27(1) (2017), 165-233.

\bibitem{Michelsohn} M.L. Michelsohn, {\it On the existence of special metrics in complex geometry}, Acta Math.149, (1982) 261-295.

\bibitem{Mue} R. M\"uller,
{\it Ricci flow coupled with harmonic map heat flow},
PhD Dissertation, 2009.

\bibitem{MS1} A. Murcia and C.S. Shahbazi,
{\it Parallel spinor flows on three-dimensional Cauchy hypersurfaces}, arXiv: 2109.13906.

\bibitem{MS2} A. Murcia and C.S. Shahbazi, {\it Parallel spinors on globally hyperbolic Lorentzian four-manifolds}, arXiv: 2011.02423, to appear in Math. Z.

\bibitem{PPZ1}
D.H Phong, S. Picard, and X.W. Zhang, {\it Geometric flows and Strominger systems}, Math. Z. 288 (2018), 101-113.

\bibitem{PPZ2}
D.H Phong, S. Picard, and X.W. Zhang, {\it Anomaly flows}, Comm. Anal. Geom. 26 (2018) No. 4, 955-1008.

\bibitem{PPZ3}
D.H. Phong, S. Picard, and X.W. Zhang, {\it The Anomaly flow and the Fu-Yau equation}, Ann. PDE 4 (2018) No. 2, Art. 13, 60 pp.

\bibitem{PPZ4}
D.H. Phong, S. Picard, and X.W. Zhang, {\it A flow of conformally balanced metrics with K\"ahler fixed points}, Math. Ann. 374 (2019) no. 3-4, 2005-2040.

\bibitem{PPZ5} D.H. Phong, S. Picard, and X.W. Zhang,
{\it The Anomaly flow on unimodular Lie groups}, Advances in Complex Geometry, 217-237, Contemp. Math. 735, Amer. Math. Soc., Providence, RI, 2019.

\bibitem{PSS}
D.H. Phong, N. Sesum, and J. Sturm,
{\it Multiplier ideal sheaves and the K\"ahler-Ricci flow},
Comm. Analysis Geom. 15 (2007) 613-632.

\bibitem{PS} S. Picard and C. Suan,
{\it Flows of $G_2$ structures associated to Calabi-Yau manifolds},
arXiv: 2209.03411.

\bibitem{PZ} S. Picard and X.W. Zhang,
{\it Parabolic Monge-Amp\`ere equations on compact K\"ahler manifolds}, 
arXiv: 1906.10235.

\bibitem{PvN} C. Pope and P. van Nieuwenhuizen, {\it Compactifications of d = 11 supergravity on K¨ahler
manifolds}, Commun. Math. Phys. 122 (1989) 281-292.

\bibitem{PW} C. Pope and N. Warner, {\it An SU(4) invariant compactification of d = 11 supergravity on a
stretched seven-sphere}, Phys. Lett. 150 B no 5, 352-356.

\bibitem{PU} M. Pujia and L. Ugarte,
{\it The Anomaly flow on nilmanifolds},
arXiv: 2004.06744.

\bibitem{Raf} A. Raffero, {\it Special solutions for the Type IIA flow},
arXiv: 2107.12264.

\bibitem{KJS} K. Smith,
Columbia University PhD thesis, 2023.

\bibitem{SoT} J. Song and G. Tian,
{\it Canonical measures and K\"ahler-Ricci flow},
J. Amer. Math. Soc. 25 (2) (2012) 303-353.

\bibitem{ST} J. Streets and G. Tian, {\it Hermitian curvature flow}, J. Eur. Math. Soc. 13 (2011), 601-634.

\bibitem{S} A. Strominger, {\it Superstrings with torsion}, Nuclear Phys. B 274 (1986), no. 2, 253-284.

\bibitem{Sze} G. Szekelyhidi,
{\it Degenerations of $C^n$ and Calabi-Yau metrics}, arXiv: 1706.00357.

\bibitem{Ta} Y. Tanaka, {\it Some boundedness properties of solutions to the Vafa-Witten equations on closed 4-manifolds}, arXiv: 1308.0862.

\bibitem{Taubes} C.H. Taubes, {\it On the singular sets of solutions to the Kapustin-Witten equations and the Vafa-Witten equations on compact K\"ahler surfaces}, arXiv: 1510.07739.

\bibitem{T} A. Tomasiello, 
{\it Reformulating supersymmetry with a generalized Dolbeault operator}, J. High
Energy Phys. (2008) no. 2, 010, arXiv:0704.2613.

\bibitem{TV} A. Tomassini and L. Vezzoni,
{\it On symplectic half-flat manifolds}, Manuscripta Math. 125 (2008) no. 4, 515-530.

\bibitem{To} P. Townsend, {\it The eleven-dimensional supermembrane revisited}, Phys. Lett. B 350 (1995), 184 – 188.

\bibitem{TY1} L.S. Tseng and S.T. Yau,
{\it Generalized cohomologies and supersymmetry}, Comm. Math. Phys. 326 (2014), no. 3, 875–885.

\bibitem{TY2}
L.S. Tseng and S.-T. Yau, {\it Cohomology and Hodge theory on symplectic manifolds}: I, J. Differential Geom. 91 (2012), no. 3, 383–416.

\bibitem{TY3} 
L.S. Tseng and S.T. Yau, {\it Cohomology and Hodge theory on symplectic manifolds}: II, J. Differential Geom. 91 (2012), no. 3, 417–443.

\bibitem{UY} K. Uhlenbeck and S.T. Yau, {\it On the existence of Hermitian-Yang-Mills connections in stable vector bundles}, Comm. Pure Appl. Math. 39 (1986), no. S, suppl., S257-S293. Frontiers of the mathematical sciences: 1985 (New York, 1985).

\bibitem{U} Y. Ustinovskyi, 
{\it Hermitian curvature flow and curvature positivity conditions}, Princeton University PhD Thesis, June 2018.

\bibitem{VW} C. Vafa and E. Witten,
{\it A strong coupling test of S-duality}, Nucl. Phys. B 431 (1994) 3-77.

\bibitem{W} E. Witten, {\it String theory dynamics in various dimensions}, Nucl. Phys. B 443 (1995), 85-126.

\bibitem{Y} S.T. Yau, {\it On the Ricci curvature of a compact K\"ahler manifold and the complex Monge-Amp\`ere equation.} I, Comm. Pure Appl. Math. 31 (1978) 339-411.



}

\end{thebibliography}
\end{document}